\documentclass[10pt]{article}
\usepackage{latexsym}
\usepackage{amsfonts}
\usepackage{enumerate}
\usepackage{multicol}
\usepackage{graphicx}
\usepackage{amssymb}
\usepackage{amsmath}
\usepackage{epic}

\title{A Walk Through Some Newer Parts of Additive Combinatorics}

\author{B\'{e}la Bajnok  \\
{\small Department of Mathematics, Gettysburg College, U.S.A.}  \\ {\small Email: bbajnok@gettysburg.edu.}}

\date{May 11, 2022}

\newtheorem{thm}{Theorem}

\newtheorem{conj}[thm]{Conjecture}

\newtheorem{prob}[thm]{Problem}

\begin{document}

\maketitle

\begin{abstract} [This paper is based on the second of two plenary talks given by the author at the 53rd Southeastern International Conference on Combinatorics, Graph Theory \& Computing, held at Florida Atlantic University,  March 7--11, 2022.]  In this survey paper we discuss some recent results and related open questions in additive combinatorics, in particular, questions about sumsets in finite abelian groups.

\end{abstract}

\thispagestyle{empty}

\section{Introduction}

We embark on a tour through some newer parts of additive combinatorics, visiting some recent results and related open questions.  Our context will be within finite abelian groups, written additively.  When our group is cyclic and of order $n$, we identify it with $\mathbb{Z}_n=\mathbb{Z}/n\mathbb{Z}$; we consider $0,1,\dots,n-1$ interchangeably as integers and as elements of $\mathbb{Z}_n$.

Much of additive combinatorics can be described as the study of combinatorial properties of sumsets.  Given a nonnegative integer $h$ and an $m$-subset $A=\{a_1, \dots, a_m\}$ of $G$, we recall the following definitions and notations:

\begin{itemize}
\item the $h$-fold {\em sumset} of $A$:
$$h A = \left \{\lambda_1 \cdot a_1 + \cdots + \lambda_m \cdot a_m \mid \lambda_i \in \mathbb{N}_0, \; \Sigma_{i=1}^m  \lambda_i  = h \right\};$$

\item the $h$-fold {\em restricted sumset} of $A$:
$$h \hat{\;} A = \left \{\lambda_1 \cdot a_1 + \cdots + \lambda_m \cdot a_m \mid \lambda_i \in \{0,1\}, \; \Sigma_{i=1}^m  \lambda_i  = h \right\};$$

\item the $h$-fold {\em signed sumset} of $A$: 
$$h_{\pm} A= \left \{\lambda_1 \cdot a_1 + \cdots + \lambda_m \cdot a_m \mid \lambda_i \in \mathbb{Z}, \; \Sigma_{i=1}^m | \lambda_i | = h \right\};$$

\item the $h$-fold {\em restricted signed sumset} of $A$: 
$$h \hat{_\pm} A = \left \{\lambda_1 \cdot a_1 + \cdots + \lambda_m \cdot a_m \mid \lambda_i \in \{0, \pm 1\}, \; \Sigma_{i=1}^m | \lambda_i | = h \right\} .$$
\end{itemize}
These four types of sumsets can be illustrated in the following diagram:

\begin{center}
\begin{tabular}{l|ccc}
& Repetition of terms & & Terms must be \\ 
& is allowed & & distinct \\ \hline \\
Terms can be & & & \\
added only & \; $hA$ & $\supseteq$ & \; $h\hat{\;}A$ \\ \\ 
& $| \bigcap$ & & $| \bigcap$  \\
Terms can be & & & \\ 
added or subtracted & \; $h_{\pm}A$  & $\supseteq$ & \; $h \hat{_{\pm}}A$ \\ \\
\end{tabular}
\end{center}
The sizes of these sumsets may vary greatly: for example, the set $A=\{2,3,5,7\}$ in $\mathbb{Z}_{53}$ (in recognition of the 53rd conference) has 
$|3A| = 14$,  $|3 \hat{\;} A|=4$, $|3_{\pm}A|=39$,  and $|3 \hat{_{\pm}}A| = 23$.

In the four sections below we discuss our favorite recent open questions about minimum sumset size, perfect bases and spanning, $(k, \ell)$-sumfree sets, and maximum-size nonbases, respectively.

For additional results and open questions on these and other related topics in additive combinatorics, we recommend the author's book \cite{Baj:2018a}.


\section{Minimum sumset size}

In this section we aim to address the following general questions: Among the subsets of a finite abelian group, all of a same given size, what is the smallest possible size that their sumsets can be?  And, conversely, how can one characterize the subsets of the group that achieve this minimum size?

Specifically, given a finite abelian group $G$ and positive integers $m$ and $h$, we introduce the following notations: 
$$\rho(G, m, h) = \min \{ |hA| \; : \; A \subseteq G, |A|=m\};$$
functions $\rho \hat{\;} (G, m, h)$, $\rho_{\pm}(G, m, h), $ and $\rho \hat{_{\pm}}(G, m, h)$ are defined analogously.  

Among these four quantities, $\rho(G, m, h)$ has the longest history and it is the only one that is fully known at the present time.  The story starts with Cauchy's result from more than two hundred years ago, which was rediscovered by Davenport a century later, so is now called the Cauchy--Davenport Theorem; we state this result using our notation, as follows.

\begin{thm}[Cauchy, cf.~\cite{Cau:1813a}; Davenport, cf.~\cite{Dav:1935a, Dav:1947a}]

For any prime $p$ and positive integer $m$, we have $$\rho(\mathbb{Z}_p, m, 2)=\min \{p, 2m-1\}.$$

\end{thm}
The general value of $\rho(G, m, h)$, which is what we explain next, has only been known for about a decade and a half.   

To start the discussion, we ask: when do subsets $A$ of $G$ have small sumsets?  There are two ideas that come to mind.
\begin{itemize}
\item Put $A$ in a coset of a subgroup:  Indeed, if $A$ is a subset of $a+H$ for some subgroup $H$ of $G$ and element $a$ of $A$, then $ hA$ is a subset of $h \cdot a+H$, and thus can have size at most $|H|$.

\item Put $A$ into an arithmetic progression:  If $A$ is a subset of $\cup_{i=0}^{k-1} \{a+i \cdot g\}$ for some $a \in A$, $g \in G$, and positive integer $k$, then  $hA$ is a subset of $ \cup_{i=0}^{hk-h} \{a+i \cdot g\}$, and so it can have size at most $hk-h+1$.

\end{itemize}
    For the case of $m$-subsets in groups of prime order $p$, these two ideas yield the upper bound $\min \{p, hm-h+1\}.$  In other groups, we may combine these ideas, and place $A$ inside an arithmetic progression of cosets of some subgroup.  

Let us carry this out more carefully.  In order to simplify notations, we assume that $A$ is an $m$-subset of the cyclic group $G=\mathbb{Z}_n$.  
First, we choose a subgroup $H$ of $G$; if $H$ is of order $d$, then we have $H=\cup_{j=0}^{d-1}\{j \cdot n/d\}$.  We will use exactly $\lceil m/d \rceil$ cosets of $H$ that form an arithmetic progression; in particular, if $m=cd+k$ for some integers $c$ and $k$ with $1 \leq k \leq d$, then we let $A$ consist of the `first' $c$ cosets of $H$, plus the `first' $k$ elements of the $(c+1)$-st coset:
$$A=\cup_{i=0}^{c-1} (i+H) \bigcup \cup_{j=0}^{k-1}\{c+j \cdot n/d\}.$$
Then $A$ has size $m=cd+k$, and we have 
$$hA=\cup_{i=0}^{hc-1} (i+H) \bigcup \cup_{j=0}^{hk-h}\{hc+j \cdot n/d\},$$ and thus $hA$ has size
\begin{eqnarray*} 
|hA| & = &  \min\{n,hcd+\min\{d,hk-h+1\}\} \\    
& = & \min\{n,hcd+d, hcd+ hk-h+1\} \\     
& = & \min\{n,(hc+1)d, hm-h+1\} \\    
& = & \min\{n,\left( h\lceil m/d \rceil -h+1 \right) \cdot d, hm-h+1\} .
\end{eqnarray*}

Introducing the notation $$f_d=\left( h\lceil m/d \rceil -h+1 \right) \cdot d,$$ we can see that 
$f_n=  n$ and $f_1=  hm-h+1$, so 
$$|hA|=\min\{f_n, f_d, f_1\}.$$
Therefore,
$$\rho (\mathbb{Z}_n, m, h) \leq \min\{\left( h\lceil m/d \rceil -h+1 \right) \cdot d \; : \; d|n \}.$$ Using Kneser's Theorem (see \cite{Kne:1953a} or \cite{Nat:1996a}), one can prove that we cannot do better, and we have the following result.

\begin{thm}[Plagne, cf.~\cite{Pla:2006a}] \label{FAU Plagne rho}

For all finite abelian groups $G$ of order $n$ and integers $m$ and $h$, we have $$\rho (G, m, h)=\min\{\left( h\lceil m/d \rceil -h+1 \right) \cdot d \; : \; d|n \}.$$

\end{thm}

Having answered the direct problem of finding $\rho (G, m, h)$, we now turn to the inverse question: what can we say about $m$-subsets $A$ of $G$ whose $h$-fold sumset has size $\rho (G, m, h)$?  The answer in groups of prime order is a consequence of Vosper's Theorem, and can be stated as follows.  

\begin{thm} [Vosper, cf.~\cite{Vos:1956a, Vos:1956b}] \label{FAU Vosper}
Suppose that $G$ is of prime order $p$.  Let $m$ and $h$ be positive integers so that $p > hm-h+1$, and suppose that $A$ is an $m$-subset of $G$.  Then $hA$ has size $hm-h+1$ if, and only if, $h=1$ or $A$ is an arithmetic progression.    
\end{thm}

The characterization in other groups is likely to be hard in its most general form.  For instance, while the only type of $6$-subset of $\mathbb{Z}_{15}$ with a 2-fold sumset of size $\rho( \mathbb{Z}_{15},6,2)=9$ consists of two cosets of the subgroup of order 3, 7-subsets with 2-fold sumsets of size $\rho( \mathbb{Z}_{15},7,2)=13$ may come in a variety of forms: unions of two cosets of the subgroup of order 3 plus one element, a coset of the subgroup of order 5 plus two other elements, or an arithmetic progression of size 7.  

As a modest generalization of Theorem \ref{FAU Vosper}, we consider here the case when $p$ is the smallest prime divisor of the order of $G$, and $m \leq p$.  Note that, in this case, from Theorem \ref{FAU Plagne rho} we have 
$$\rho (G, m, h)=\min\{p, hm-h+1\}.$$
It is not hard to see that the following results follow from Kemperman's Theorem (cf.~\cite{Kem:1960a}) and Kneser's Theorem, respectively.

\begin{thm}  \label{FAU Kemperman}
Let $p$ be the smallest prime divisor of the order of $G$.  Let $m$ and $h$ be positive integers so that $p > hm-h+1$, and suppose that $A$ is an $m$-subset of $G$.  Then $hA$ has size $hm-h+1$ if, and only if, $h=1$ or $A$ is an arithmetic progression.
\end{thm}

\begin{thm} \label{FAU follows from Kneser}
Let $p$ be the smallest prime divisor of the order of $G$.  Let $m$ and $h$ be positive integers so that $m \leq p < hm-h+1$, and suppose that $A$ is an $m$-subset of $G$.  Then $hA$ has size $p$ if, and only if, $A$ is contained in a coset of some subgroup $H$ of $G$ with $|H|=p$.   
\end{thm}

The general inverse question remains largely open.

\begin{prob}
For each abelian group $G$ and positive integers $m$ and $h$, find a characterization of $m$-subsets of $G$ that have $h$-fold sumsets of size $\rho (G, m ,h)$.

\end{prob}

We now turn to minimum sumset size for restricted addition.  The value of
$$\rho \hat{\;} (G, m, h) = \min \{ |h \hat{\;} A| \; : \; A \subseteq G, |A|=m\}$$
is largely unknown.  However, after several decades of being known as the Erd\H{o}s--Heilbronn Conjecture, the value is finally known in groups of prime order. 

\begin{thm}[Dias da Silva and Hamidoune, cf.~\cite{DiaHam:1994a}; Alon, Nathanson, and Ruzsa, cf.~\cite{AloNatRuz:1995a, AloNatRuz:1996a}]

For all primes $p$ and positive integers $m$ and $h$, we have $$ \rho \hat{\;} (\mathbb{Z}_p, m, h) = \min\{p,hm-h^2+1\}.$$

\end{thm}
  
Regarding groups of composite order, we can provide the following sharp result.

\begin{thm} [\cite{Baj:2013a}] \label{FAU rho hat 2}

For all positive integers $n$, $m$, and $h$ we have 
$$ \rho\hat{\;} (\mathbb{Z}_n,m,2) \le \left\{
\begin{array}{ll}
\min\{\rho(\mathbb{Z}_n,m,2), 2m-4\} & \mbox{if} \; 2|n \; \mbox{and} \;  2|m, \mbox{or} \; (2m-2)|n \; \mbox{and} \;  m-1 \neq 2^k; \\ \\
\min\{\rho(\mathbb{Z}_n,m,2), 2m-3\} & \mbox{otherwise.} 
\end{array}
\right.$$

\end{thm}

We believe that $ \rho\hat{\;} (\mathbb{Z}_n,m,2)$ is given exactly in Theorem \ref{FAU rho hat 2}, so we ask the following.

\begin{prob}  \label{FAU rho hat 2 prob}

Prove that equality holds in Theorem \ref{FAU rho hat 2}.
\end{prob} 

Regarding general abelian groups, we mention the following conjecture.

\begin{conj} [Lev, cf.~\cite{Lev:2000a}]  \label{FAU Lev conj}
For any abelian group $G$ and for all positive integers $m$, we have $$\rho \hat{\;} (G, m, 2) \geq \min \{ \rho (G, m, 2), 2m-2-|L|\}$$ 
where $L$ is the subgroup of involutions in $G$.  
\end{conj}
We may observe that in a cyclic group of order $n$ we have $|L|=2$ when $n$ is even and $|L|=1$ when $n$ is  odd, cf. Theorem \ref{FAU rho hat 2} above.  We also mention here that Conjecture \ref{FAU Lev conj} was established for elementary abelian groups:

\begin{thm} [Eliahou and Kervaire, cf.~\cite{EliKer:1998a}]

For all positive integers $m$ and $r$ and for odd primes $p$, we have $$\rho \hat{\;} (\mathbb{Z}_p^r, m, 2) \geq \min \{ \rho (\mathbb{Z}_p^r, m, 2), 2m-3\}.$$ 

\end{thm}
Note that since $2A=2 \hat{\;}A \cup \{0\}$ for any subset $A$ of $\mathbb{Z}_2^r$, we trivially have 
$$\rho \hat{\;} (\mathbb{Z}_2^r, m, 2)= \rho(\mathbb{Z}_2^r, m, 2) -1.$$

The general problem of finding $\rho \hat{\;} (G, m, h)$ is largely open and likely to be difficult: as we see, unlike $\rho (G, m, h)$, it depends on the structure of $G$ and not just the order of $G$.  As a potential next case, we have the following.

\begin{conj} \label{FAU conj p smallest}
Let $p$ be the smallest prime divisor of the order of $G$.  If $m$ and $h$ are positive integers so that $m \le p$, then $$\rho \hat{\;} (G, m, h)=\min\{p, hm-h^2+1\}.$$
\end{conj}

The inverse problems paralleling Theorems \ref{FAU Kemperman} and \ref{FAU follows from Kneser} above are both conjectures for restricted addition.

\begin{conj} \label{FAU hat conj}
Let $p$ be the smallest prime divisor of the order of $G$.  Let $m$ and $h$ be positive integers so that $p > hm-h^2+1$, and suppose that $A$ is an $m$-subset of $G$.  Then $hA$ has size $hm-h^2+1$ if, and only if, $h=1$; $A$ is an arithmetic progression; or $h=2$, $m=4$, and $A=\{a,a+g_1,a+g_2,a+g_1+g_2\}$ for some $a \in A$ and $g_1, g_2 \in G$.
\end{conj}

\begin{conj} 
Let $p$ be the smallest prime divisor of the order of $G$.  Let $m$ and $h$ be positive integers so that $m \leq p < hm-h^2+1$, and suppose that $A$ is an $m$-subset of $G$.  Then $hA$ has size $p$ if, and only if, $A$ is contained in a coset of some subgroup $H$ of $G$ with $|H|=p$.   
\end{conj}

K\'arolyi proved that Conjectures \ref{FAU conj p smallest} and \ref{FAU hat conj} hold for $h=2$, cf.~\cite{Kar:2003a, Kar:2004a, Kar:2005a}.

Let us now proceed to signed sumsets and the question of finding 
$$\rho_{\pm}(G, m, h) = \min \{ |h_{\pm}A|  \; : \; A \subseteq G, |A|=m\}.$$
The value of this function is not generally known, though we have some interesting partial results, and even a conjecture for all cases.  Like we have just seen with $\rho \hat{\;} (G, m, h)$, the value of $\rho_{\pm} (G, m, h)$ depends on the structure of $G$, not just the order of $G$.  We should mention that, perhaps surprisingly, $\rho_{\pm} (G, m, h)$ often agrees with $\rho (G, m, h)$, although generally $h_{\pm}A$ is a lot larger than $hA$.  
For example, among groups of order 24 or less, we have $\rho_{\pm} (G, m, h)=\rho (G, m, h)$, with the only  exception of $G=\mathbb{Z}_3^2$, $m=4$, and $h=2$, where $\rho_{\pm} (\mathbb{Z}_3^2, 4,2)=8$ but $\rho (\mathbb{Z}_3^2, 4,2)=7$.  The two functions always agree for cyclic groups, though:

\begin{thm} [with Matzke, cf.~\cite{BajMat:2014a}] \label{FAU Matzke cyclic}

For any cyclic group $G$ and for all positive integers $m$ and $h$, we have $$\rho_{\pm} (G, m, h)=\rho (G, m, h).$$

\end{thm}
For the proof of Theorem \ref{FAU Matzke cyclic}, we constructed a subset $R \subseteq G$ so that $R$ is symmetric (that is, $R = -R$ and thus $h_{\pm}R=hR$), has size at least $m$, and $$|h_{\pm}R| \leq \left( h\lceil m/d \rceil -h+1 \right) \cdot d.$$  The result now follows from Theorem \ref{FAU Plagne rho}.

For noncyclic groups, the situation is substantially more complicated, though \cite{BajMat:2014a} contains a conjecture for $\rho_{\pm} (G, m, h)$ in all cases.  Here we only discuss elementary abelian groups $\mathbb{Z}_p^r$; since we clearly have 
$$\rho_{\pm} (\mathbb{Z}_2^r, m, h) = \rho (\mathbb{Z}_2^r, m, h),$$
we assume that $p$ is an odd prime.  We have the following rather delicate result.

\begin{thm} [with Matzke, cf.~\cite{BajMat:2014b}] \label{FAU Matzke p group}

Suppose that $p$ is an odd prime and that $m$, $h$, and $r$ are positive integers; we use the notation $f_1=hm-h+1$.  We define $\delta$ to be 0 or 1, depending on whether $p-1$ is divisible by $h$ or not.  We then let $k$ to be the largest integer for which  $p^k  +\delta$ is at most $f_1$, and then set $c$ to be the largest integer for which $(hc+1) \cdot p^k + \delta$ is at most $f_1$.  

If $m \leq (c+1) \cdot p^k$, then $$\rho_{\pm} (\mathbb{Z}_p^r, m, h) = \rho (\mathbb{Z}_p^r, m, h).$$

\end{thm}
We believe that the condition for equality in Theorem \ref{FAU Matzke p group} is also necessary.

\begin{conj} \label{FAU Matzke conj p group}

If, using the notations of Theorem \ref{FAU Matzke p group}, $m > (c+1) \cdot p^k$, then $$\rho_{\pm} (\mathbb{Z}_p^r, m, h) > \rho (\mathbb{Z}_p^r, m, h).$$

\end{conj}
Using results of Vosper, Kemperman, and Lev, we have proved Conjecture \ref{FAU Matzke conj p group} for the case of $r=2$ and $h=2$.

The inverse problems regarding $\rho_{\pm} (G, m, h)$ are quite interesting as well: perhaps surprisingly, symmetric subsets don't always yield minimum signed sumset size.  Recall that a subset $A$ in a group is symmetric if $A$ equals $-A$ and asymmetric when $A$ and $-A$ are disjoint.  In addition, we define $A$ to be near-symmetric when it is possible to remove one element from it after which it becomes symmetric.  We have the following inverse result.

\begin{thm} [with Matzke, cf.~\cite{BajMat:2014a}] 

Suppose that $G$ is a finite abelian group and that $m$ and $h$ are positive integers.  Let ${\cal A}(G,m)$ denote the collection of $m$-subsets of $g$ that are symmetric, asymmetric, or near-symmetric.  We then have  
$$\rho_{\pm} (G,m,h)= \min \{|h_{\pm} A| \; : \; A \in {\cal A}(G,m)\}.$$

\end{thm}
We note that all three types of subsets are essential; it may be an interesting question to characterize all situations where the minimum signed sumset size is achieved by symmetric, asymmetric, or near-symmetric subsets, respectively.

We close this section by mentioning our fourth quantity,  
   $$\rho \hat{_{\pm}}(G, m, h) = \min \{ |hA| \; : \; A \subseteq G, |A|=m\},$$
though only to say that we know very little about it.  Clearly,
$$ \rho \hat{\;}(G, m, h) \leq \rho \hat{_{\pm}}(G, m, h) \leq \rho_{\pm}(G, m, h),$$
and both inequalities may be strict.  Perhaps a good place to start the evaluation is in cyclic groups and for $h=2$.

\begin{prob}
Evaluate $\rho \hat{_{\pm}}(\mathbb{Z}_n, m, 2)$ for positive integers $m$ and $n$.  

\end{prob}


\section{Perfect bases and spanning sets}

In this section we look for perfection: subsets of abelian groups that generate each element of the group uniquely.  Namely, for a finite abelian group $G$, a subset $A$ of $G$, and a positive integer $s$ we introduce the following definitions:

\begin{itemize}
\item $A$ is a {\em perfect $s$-basis} in $G$ if each element of $G$ can be written uniquely as a sum of at most $s$ elements of $A$.

\item $A$ is a {\em perfect restricted $s$-basis} in $G$ if each element of $G$ can be written uniquely as a sum of at most $s$ distinct elements of $A$.

\item $A$ is a {\em perfect $s$-spanning set} in $G$ if each element of $G$ can be written uniquely as a signed sum of at most $s$  elements of $A$.

\item $A$ is a {\em perfect restricted $s$-spanning set} in $G$ if each element of $G$ can be written uniquely as a signed sum of at most $s$ distinct elements of $A$.
\end{itemize}
As we will see, perfection does exist, though sometimes it may be difficult to find.

We start with perfect $s$-bases.  As far as we know, perfect bases have not been investigated yet before, though the concept of $s$-basis has enjoyed a rich history since it was first discussed by Erd\H{o}s and Tur\'an  in \cite{ErdTur:1941a} eighty years ago; see, for example, \cite{ErdGra:1980b}, \cite{ErdNat:1987a}, \cite{LamThaPla:2017a}, \cite{Nat:1974a},  \cite{Nat:1996a}, and \cite{Nat:2014a}.

Trivially, the set of nonzero elements is a perfect $1$-basis in any group $G$, and the 1-element set consisting of a generator of $G$ is a perfect $s$-basis in the cyclic group of order $s+1$. 
It turns out that there are no others:

\begin{thm}[with Berson and Just, cf.~\cite{BajBerJus:2022a}] \label{FAU CoHo-1}
If a subset $A$ of a finite abelian group $G$ is a perfect $s$-basis, then $s=1$ and $A=G \setminus \{0\}$, or $G$ is cyclic of order ${s+1}$ and $A$ consists of a single element.

\end{thm}

The proof of Theorem \ref{FAU CoHo-1} is based on the fact that if $A$ is a perfect $s$-basis in $G$, then the $s+1$ subsets 
$$-A, \; A-A, \; 2A-A, \; \dots, \; (s-1)A-A, \; \mbox{and} \; (s-1)A$$ are pairwise disjoint.  We can then compute that 
$$|\cup_{h=0}^{s-1} (hA-A) \; \cup \; (s-1)A| =  {m+s \choose s} + \frac{(m-1)(s-1)}{s} {m+s-2 \choose s-1}.$$  However, if $A$ is a perfect $s$-basis of size $m$ in $G$, then 
$$ |G|=\sum_{h=0}^s |hA|= \sum_{h=0}^s {m+h-1 \choose h} = {m+s \choose s},$$ and this is only possible if $m=1$ or $s=1$, as claimed. 

Moving on to perfect restricted $s$-bases, we first observe that the instances of $|A|=1$ or $s=1$ are identical for restricted addition and unrestricted addition, and are as listed in Theorem \ref{FAU CoHo-1}.  As a major contrast, however, there are infinitely many perfect restricted $s$-bases with $|A| \geq 2$ and $s \geq 2$, as we explain below.  We are able to provide a complete characterization for the case of $s=2$, for which we have the following result.

\begin{thm} [with Berson and Just, cf.~\cite{BajBerJus:2022a}] \label{FAU CoHo-3}
A finite abelian group $G$ has a perfect restricted $2$-basis if, and only if, it is isomorphic to one of the following groups: $\mathbb{Z}_2, \; \mathbb{Z}_4, \; \mathbb{Z}_7, \; \mathbb{Z}_2^2, \; \mathbb{Z}_2^4$, or $\mathbb{Z}_2^2 \times \mathbb{Z}_4$.

\end{thm}

Our strategy for proving Theorem \ref{FAU CoHo-3} is similar to the unrestricted case, but exhibiting sets that are pairwise disjoint seems more elusive, hence our result is only for $s =2$.  For a given perfect restricted $2$-basis $A$, we consider the set
$$T=(A-A) \cup A,$$ and prove that, unless $G$ isomorphic to $\mathbb{Z}_4$, $\mathbb{Z}_7$, or $\mathbb{Z}_2^2 \times \mathbb{Z}_4$ (in which cases perfect restricted $2$-bases exist) or to an elementary abelian 2-group, $T$ has size more than the order of $G$.      
However, if all nonzero elements of $G$ have order 2, then $$A-A=2A=\{0\} \cup 2 \hat{\;} A,$$ and thus $T=\cup_{h=0}^2 h\hat{\;} A.$ Therefore, $A$ being a perfect restricted $2$-basis in $G$ is equivalent to having $T=G$, resulting in no contradiction.  

Luckily, a problem of Ramanujan comes to the rescue. In 1913, Ramanujan asked in \cite{Ram:2000a} whether the quantity $2^k-7$ can be a square number for any integer $k$ besides 3, 4, 5, 7, and 15 (see also Question 464 in \cite{BerChoKan:1999a}).   The negative answer was given by Nagell in 1948 (see \cite{Nag:1948a}; \cite{Nag:1961a} for the English version).  Suppose now that $G$ is the elementary abelian 2-group of rank $r$ that has a perfect restricted $2$-basis of size $m$: we then must have 
$$2^r={m \choose 0} + {m \choose 1} + {m \choose 2},$$ and therefore $$2^{r+3}-7=(2m+1)^2.$$  We thus have exactly four choices for $r$: 1, 2, 4, and $12$.  Perfect restricted $2$-bases exist for the first three choices, and we were able to show that they do not in $\mathbb{Z}_2^{12}$, completing the proof of Theorem \ref{FAU CoHo-3}.  

We know of perfect restricted $s$-bases of size $m$ exist for higher $s$ as well: 

\begin{itemize}
  \item for $m \leq s$, in every group of order $|G|=2^m$;
  \item for $m=s+1$, in  $G \cong \mathbb{Z}_{2^{s+1}-1}$; and
  \item for $m=2s+1$, in $G \cong \mathbb{Z}_2^{2s}$, as well as in $G \cong \mathbb{Z}^{2s-2} \times \mathbb{Z}_4$.
\end{itemize}
We believe that there are no others.

\begin{conj}
For integers $m \geq 2$ and $s \geq 2$, perfect restricted $s$-bases of size $m$ exist if, and only if, $G$ has order $2^m$ or is isomorphic to $\mathbb{Z}_{2^{s+1}-1}$, $\mathbb{Z}_2^{2s}$, or to $\mathbb{Z}^{2s-2} \times \mathbb{Z}_4$.

\end{conj}

Next, we turn to perfect spanning sets, that is, subsets $A$ of our group $G$ where each element of $G$ can be written uniquely as a signed sum of at most $s$  elements of $A$.  Note that, if $A$ is a perfect $s$-spanning set of $G$ and is of size $m$, then the order of $G$ must equal the Delannoy number $$a(m,s)=\sum_{i\geq 0} {s \choose i} \cdot {m \choose i} \cdot 2^i.$$  
(Delannoy numbers may also be defined by the recursion
$$a(m,s) = a(m-1,s)+a(m-1,s-1)+a(m,s-1),$$ together with the initial conditions of $a(m,0)=a(0,s)=1.$)
We are aware of only the following perfect spanning sets:

\begin{itemize}
  \item 
$s=1$ (and $m$ arbitrary), in which case $|G|= a(m,1)=2m+1$:  $A$ is perfect in $G$ if, and only if, $A$ and $-A$ partition $G \setminus \{0\}$.

\item

$m=1$ (and $s$ arbitrary), in which case  $|G|=a(1,s)=2s+1$: 
$A$ is perfect in $G$ if, and only if,  $G$ is cyclic of order $2s+1$, and 
$A=\{a\}$ with $\gcd (a,|G|)=1$.

\item

$m=2$ (and $s$ arbitrary), in which case  
$|G|=a(2,s)=2s^2+2s+1$: 
$A$ is perfect in $G$ if $G$ is cyclic of order $2s^2+2s+1$ and
 $A=c \cdot \{s,s+1\}$ with $\gcd (c,|G|)=1$.
\end{itemize}
We believe that this rather short list is complete:

\begin{conj} \label{FAU conj indep span}

The only instances of perfect sets in $G$ are those three just mentioned.  In particular, we must have $s=1$ or $m \in \{1,2\}$.

\end{conj}

At the present time, we know little about perfect restricted spanning sets, though we suspect that their theory (with the powers of 3 replacing the powers of 2) can be developed analogously to perfect restricted bases.


\section{$(k,l)$-sumfree sets}

Recall that a subset $A$ of $G$ is called {\em sumfree} if it does not contain the sum of two (not necessarily distinct) of its elements, that is, if $A$ and $2A$ are disjoint.  More generally, for positive integers $k$ and $l$, with $k>l$, we call a subset $A$ of $G$ {\em $(k,l)$-sumfree} if $kA$ and $lA$ are disjoint.  The main question we intend to discuss in this section is how large a $(k,l)$-sumfree subset of a given group can be.

\begin{prob}
For each pair of positive integers $k$ and $l$, and for each finite abelian group $G$, find the maximum size $\mu (G, \{k,l\})$ of $(k,l)$-sumfree subsets in $G$.

\end{prob}

Sumfree sets in abelian groups were first introduced by Erd\H{o}s in \cite{Erd:1965a} and then studied systematically by Wallis, Street, and Wallis in \cite{WalStrWal:1972a}.   

We can construct sumfree sets in $G$ by selecting a subgroup $H$ in $G$ for which $G/H$ is cyclic and then taking the `middle one-third' of the cosets of $H$.  More precisely, with $d$ denoting the index of $H$ in $G$, we see that 
$$\mu (G, \{2,1\}) \ge \max_{d |e(G)} \left \{ \left \lceil \frac{d-1}{3} \right \rceil \cdot \frac{n}{d}  \right \},$$ where $e(G)$ is the exponent of $G$.  
Using a version of Kneser's Theorem, Diamanda and Yap (see \cite{DiaYap:1969a} or \cite{WalStrWal:1972a}) proved in the 1960s that we cannot do better in cyclic groups; this result was only extended to the general case in 2005:   
\begin{thm} [Green and Ruzsa, cf.~\cite{GreRuz:2005a}]  \label{FAU Green and Ruzsa thm}
For any abelian group $G$ of order $n$ and exponent $e(G)$, we have 
$$
\mu (G, \{2,1\}) = \max_{d |e(G)} \left \{ \left \lceil \frac{d-1}{3} \right \rceil \cdot \frac{n}{d}  \right \}.
$$
\end{thm}
We should note that the proof of this result relies in part on a computer program.

For general $k$ and $l$, the first result was given by Bier and Chin (cf.~\cite{BieChi:2001a}) who evaluated $\mu (\mathbb{Z}_p, \{k,l\})$ for prime values of $p$.  
This was generalized by Hamidoune and Plagne for general cyclic groups, but only when their order was relatively prime to $k-l$:

\begin{thm} [Hamidoune and Plagne, cf.~\cite{HamPla:2003a}]  \label{FAU mu when n and k-l are rel prime}
If $k-l$ is relatively prime to $n$, then 
$$
\mu (\mathbb{Z}_n, \{k,l\}) = \max_{d |n} \left \{ \left \lceil \frac{d-1}{k+l} \right \rceil \cdot \frac{n}{d}  \right \}.
$$
\end{thm}

Just recently, we were able to generalize Theorem \ref{FAU mu when n and k-l are rel prime} to the case when $n$ and $k-l$ are not relatively prime.

\begin{thm} [with Matzke, cf.~\cite{BajMat:2019a}] \label{FAU mu all n k l}
For all positive integers $n$, $k$, and $l$, with $k>l$, we have
$$ \mu (\mathbb{Z}_n, \{k,l\}) = \max_{d |n} \left \{ \left \lceil \frac{d-(\delta-r)}{k+l} \right \rceil \cdot \frac{n}{d}  \right \},
$$
where $\delta=\gcd(d,k-l)$ and $r$ is the remainder of $l  \lceil (d-\delta)/(k+l) \rceil$ mod $\delta$.
\end{thm}

Let us review our approach for the proof of Theorem \ref{FAU mu all n k l}.  The main role is played by arithmetic progressions, that is, sets of the form
$$A=\{a+ i \cdot b \mid i=0,1,\dots,m-1\}$$ for some positive integer $m$ and elements $a$ and $b$ of $\mathbb{Z}_n$. In \cite{HamPla:2003a}, Hamidoune and Plagne proved that, if $n$ and $k-l$ are relatively prime, then $\mu (\mathbb{Z}_n, \{k,l\})$ equals     
$$\max_{d |n} \left\{ \alpha (\mathbb{Z}_d, \{k,l\}) \cdot \frac{n}{d}  \right \},$$ where $\alpha (\mathbb{Z}_d, \{k,l\})$ is the maximum size of a $(k,l)$-sumfree arithmetic progression in $\mathbb{Z}_d$.  Although the situation is considerably more intricate when $n$ and $k-l$ are not relatively prime, it turns out (see \cite{Baj:2009a}) that the identity remains valid.  This then reduces the problem of finding $ \mu (\mathbb{Z}_n, \{k,l\}) $ to arithmetic progressions only.

When trying to evaluate $\alpha (\mathbb{Z}_d, \{k,l\})$, one naturally considers two types of arithmetic progressions: those with a common difference that is not relatively prime to $d$ (in which case the set is contained in a coset of a proper subgroup), and those where the common difference is relatively prime to $d$ (in which case the set, unless of size $1$, is not contained in a coset of a proper subgroup).  The main difficulty is caused by those of the first type; luckily, however, we were able to prove that it does not matter:

\begin{thm} [with Matzke, cf.~\cite{BajMat:2019a}]  \label{FAU gamma is enough}

For all positive integers $n$, $k$, and $l$ with $k>l$ we have
$$ \mu (\mathbb{Z}_n, \{k,l\}) =  \max_{d |n} \left\{ \gamma (\mathbb{Z}_d, \{k,l\}) \cdot \frac{n}{d}  \right \}, 
$$ where $\gamma (\mathbb{Z}_d, \{k,l\})$ is the size of the largest $(k,l)$-sumfree interval in $\mathbb{Z}_d$.

\end{thm}

Since intervals are easy to handle, finding $\gamma (\mathbb{Z}_d, \{k,l\})$ is a routine calculation; this then yields Theorem \ref{FAU mu all n k l} above.

Turning now to the inverse problem of classifying all $(k,l)$-sumfree subsets of $G$ of maximum size, we first mention that the problem for $(k,l)=(2,1)$ has finally been completed in 2016 in a complicated paper by  Balasubramanian, G. Prakash, and D. S. Ramana (see \cite{BalPraRam:2016a}).  For other $k$ and $l$, we have the following result.

\begin{thm} [Plagne, cf.~\cite{Pla:2002a}]
Let $p$ be a prime and $k$ and $l$ be positive integers with $k>l$ and $k \geq 3$; assume also that $p$ does not divide $k-l$.  If $A$ is a $(k,l)$-sum-free set of size 
$$|A|=\mu (\mathbb{Z}_p, \{k,l\})=\left \lfloor \frac{p-2}{k+l} \right \rfloor +1$$ in $\mathbb{Z}_p$, then $A$ is an arithmetic progression.

\end{thm}
 The result is interesting in that it is simpler than the sumfree case (where sumfree subsets that are not arithmetic progressions may also have maximum size $\mu (\mathbb{Z}_p, \{2,1\})$; see \cite{Baj:2018a}).

\begin{prob}
For all positive integers $n$, $k$, and $l$, classify all $(k,l)$-sumfree subsets of $\mathbb{Z}_n$ that have maximum size $\mu (\mathbb{Z}_n, \{k,l\})$.

\end{prob}


\section{Sumset size of maximum-size nonbases}

Recall that a nonempty subset $A$ of $G$ is {\em $h$-complete} (alternatively, a {\em basis of order $h$}) if $hA=G$ and that if $hA$ is a proper subset of $G$, we say that $A$ is {\em $h$-incomplete}.  The {\em $h$-critical number}  $\chi  (G, h) $ of $G$ is defined as the smallest positive integer $m$ for which all $m$-subsets of $G$ are $h$-complete; that is:
$$\chi  (G, h)  =  \min \{m \; : \; A \subseteq G,  |A| \geq m \Rightarrow h  A=G \}.$$
It is easy to see that for all $G$ and $h$ we have $hG=G$, so $\chi  (G, h)$ is well defined.
The value of  ${\chi}  (G, h)$ is now known for every $G$ and $h$.

\begin{thm}[\cite{Baj:2014a}] \label{FAU critical Bajnok}

For all abelian groups $G$ of order $n$ and all positive integers $h$, we have 
$$\chi(G, h) = \max \left\{ \left( \left \lfloor \frac{d-2}{h} \right \rfloor +1  \right) \cdot \frac{n}{d}  \; : \; d|n \right\} + 1.$$

\end{thm}

In particular, by Theorem \ref{FAU critical Bajnok}, we have the following explicit formulas:
\begin{eqnarray*}
\chi (G, 1)  & = & n; \\ \\
\chi (G, 2) & = & \left\lfloor \frac{n}{2} \right\rfloor +1;  \\ \\ 
\chi (G, 3) & = & \left\{
\begin{array}{ll}
\left(1+\frac{1}{p}\right) \frac{n}{3} & \mbox{if $n$ has prime divisors congruent to 2 mod 3,} \\ & \mbox{and $p$ is the smallest such divisor,}\\ \\
\left\lfloor \frac{n}{3} \right\rfloor & \mbox{otherwise.}\\
\end{array}\right.
\end{eqnarray*}
Explicit expressions for $\chi(G, h)$ get increasingly more complicated as $h$ increases.

The question that we here try to address is the following: What can one say about the size of $hA$ if $A$ is an $h$-incomplete subset of maximum size in $G$?  Namely, we aim to determine the set
$$S(G,h)=\{ |hA| \; : \; A \subset G, \; |A|={\chi}  (G, h)-1, \; hA \neq G\}.$$ 

Trivially, $S(G,1)=\{n-1\}$.  For $h=2$ and $h=3$ we have the following results.

\begin{thm}[with Pach, cf.~\cite{BajPac:2022a}]

Let $G$ be an abelian group of order $n$.

\begin{itemize}

\item If the exponent of $G$ is divisible by $4$, then $$S(G,2)=\left \{ n-n/d \; : \; d|n, \; 2|d \right \}.$$    
\item If the exponent of $G$ is even but not divisible by $4$, then 
$$S(G,2)=\left \{ n-n/d \; : \; d|n, \; 2|d, \; d \neq 4 \right \}.$$     
\item If $n$ is odd and $n>9$, then $S(G,2)=\{n-2,n-1\}$.
\end{itemize}

\end{thm}

\begin{thm}[with Pach, cf.~\cite{BajPac:2022a}]

Let $G$ be an abelian group of order $n$.

\begin{itemize}

\item 

If $n$ has prime divisors congruent to $2$ mod $3$ and $p$ is the smallest, then $S(G,3)=\{n-n/p\}$.    

\item

If $n$ has no prime divisors congruent to $2$ mod $3$ but $3|n$, then
$$
S(G,3) = \left \{ n-n/d \; : \; d|n, \; 3|d, \; d \neq 3 \right \}  \cup \left \{ n-2n/d  \; : \; d|n, \; 1 \leq \nu_3(d) \leq \nu_3(e(G)) \right \},
$$
where $e(G)$ is the exponent of $G$, and $\nu_3(t)$ is the highest power of 3 that divides $t$.     

\item

If all divisors of $n$ are congruent to $1$ mod $3$, but $G \not \cong \mathbb{Z}_7^r$, then $S(G,3)=\{n-3,n-1\}$.

\item

If $G \cong \mathbb{Z}_7^r$, then $S(G,3)=\{n-3\}$.

\end{itemize}
\end{thm}

The analogous question for restricted addition is considerably more challenging; in fact, we do not even know the value of the restricted critical number 
$$\chi \hat{\;} (G, h)  =  \min \{m \; : \; A \subseteq G,  |A| \geq m \Rightarrow h \hat{\;} A=G \}$$
for $h \ge 3$ in all cases.  However, for $h=2$ we have $\chi \hat{\;} (G, 2)  = \lfloor n/2 \rfloor +2$ (cf.~\cite{Baj:2014a, Baj:2017a}), and the following result for 
$$S \hat{\;} (G,2)=\{ |2 \hat{\;} A| \; : \; A \subset G, \; |A|={\chi} \hat{\;} (G, 2)-1, \; 2 \hat{\;} A \neq G\}$$ in the case of cyclic groups.

\begin{thm} [Gallardo, Grekos, et al.; cf.~\cite{GalGre:2002a}] \label{thm g G et al} 
Let $G$ be a cyclic abelian group of order $n$.

\begin{itemize}

\item If $n=2^k$ for some positive integer $k$, then $S \hat{\;} (G,2) = \{ n-1\}$.
\item In all other cases, $S \hat{\;} (G,2) =\{n-2, n-1\}$.
\end{itemize}
\end{thm}

The problem of finding critical numbers $\chi \hat{\;} (G, h)$ for other $G$ and $h$ (cf.~\cite{Baj:2016b} for a review of known cases), and the corresponding sets $S \hat{\;} (G,h)$, appear to be challenging.

We hope that our short walk in additive combinatorics was enjoyable.  These and the many other interesting sights are well worth additional visits!


\begin{thebibliography}{99}





\bibitem{AloNatRuz:1995a} N. Alon, M. B. Nathanson, and I. Ruzsa, Adding distinct congruence classes modulo a prime. {\em Amer. Math. Monthly} {\bf 102} (1995), no. 3, 250--255.


\bibitem{AloNatRuz:1996a} N. Alon, M. B. Nathanson, and I. Ruzsa, The polynomial method and restricted sums of congruence classes. {\em J. Number Theory} {\bf 56} (1996), no. 2, 404--417.



 
\bibitem{Baj:2009a} B. Bajnok, On the maximum size of a $(k,l)$-sumfree subset of an abelian group.  \emph{Int. J. Number Theory}  {\bf 5} (2009), no. 6, 953--971.
 
 
 
\bibitem{Baj:2013a} B. Bajnok, On the minimum size of restricted sumsets in cyclic groups.  \emph{Acta Math. Hungar.}  {\bf 148} (2016), no. 1, 228--256.
 
\bibitem{Baj:2014a} B. Bajnok, The $h$-critical number of finite abelian groups.  {\em Unif.  Distrib.  Theory} {\bf 10} (2015), no.2, 93--15.
 
 
\bibitem{Baj:2016b}  B. Bajnok, More on the $h$-critical numbers of finite abelian groups.  {\em Ann. Univ. Sci. Budapest. E\"otv\"os Sect. Math.}  {\bf 59} (2016), 113--122.
 
\bibitem{Baj:2017a}  B. Bajnok, Corrigendum to ``The $h$-critical number of finite abelian groups.'' {\em Unif. Distrib. Theory} {\bf 12} (2017), no. 2, 119--124.

\bibitem{Baj:2017b} B. Bajnok, Open problems about sumsets in finite abelian groups: minimum sizes and critical numbers.  To appear in {\em Combinatorial and Additive Number Theory II: CANT, New York, NY, USA, 2015 and 2016,
Springer, New York,} 2017.

\bibitem{Baj:2018a} B. Bajnok, Additive Combinatorics: A Menu of Research Problems.  {\em CRC Press, Boca Raton}, 2018, xix+390 pp. 

\bibitem{BajBerJus:2022a}  B. Bajnok, C. Berson, and H.~A. Just, On Perfect Bases in Finite Abelian Groups.  To appear in {\em Involve} (2022).


\bibitem{BajMat:2014a} B. Bajnok and R. Matzke, The minimum size of signed sumsets.  {\em Electron. J. Combin.}  {\bf 22} (2015), no. 2, Paper 2.50, 17 pp.
 
\bibitem{BajMat:2014b} B. Bajnok and R. Matzke, On the minimum size of signed sumsets in elementary abelian groups.  {\em J. Number Theory} {\bf 159} (2016), 384--401.
 


\bibitem{BajMat:2019a} B. Bajnok and R. Matzke, On the maximum size of $(k,l)$-sumfree sets in cyclic groups.  {\em Bulletin of the Australian Mathematical Society} {\bf 99} (2019), no.~2, 184--194.



\bibitem{BajPac:2022a} B. Bajnok and P.~P. Pach, On sumsets of nonbases of maximum size.  Submitted (2022). 




\bibitem{BajRuz:2003a} B. Bajnok and I. Ruzsa, The independence number of a subset of an abelian group.  \emph{Integers} {\bf 3} (2003), A2, 23 pp.

  
\bibitem{BalPraRam:2016a} R. Balasubramanian, G. Prakash, and D. S. Ramana, Sum-free subsets of finite abelian groups of type III.  {\em European J. Combin.} {\bf 58} (2016), 181--202.

\bibitem{BerChoKan:1999a} B. C. Berndt,  Y-S. Choi, and S-I. Kang, The problems submitted by Ramanujan to the Journal of the Indian Mathematical Society. {\em Contemp. Math.} {\bf 236} (1999), 15--56.



\bibitem{BieChi:2001a} T. Bier and A. Y. M. Chin, On $(k,l)$-sets in cyclic groups of odd prime order. {\em Bull. Austral. Math. Soc.} {\bf 63} (2001), no. 1, 115--121.



\bibitem{Cau:1813a} A--L. Cauchy, Recherches sur les nombres. {\em J. \'Ecole Polytechnique} {\bf 9} (1813), 99--123.
 
 

 
\bibitem{Dav:1935a} H. Davenport, On the addition of residue classes. {\em J. London Math. Soc.} {\bf 10} (1935), 30--32.
 
\bibitem{Dav:1947a} H. Davenport, A historical note. {\em J. London Math. Soc.} {\bf 22} (1947), 100--101.
 
 
\bibitem{DiaYap:1969a} P. H. Diananda and H. P. Yap, Maximal sumfree sets of elements of finite groups.  {\em Proc. Japan Acad.} {\bf 45} (1969), 1--5.
 
\bibitem{DiaHam:1994a} J. A. Dias Da Silva and Y. O. Hamidoune, Cyclic space for Grassmann derivatives and additive theory.  {\em Bull. London Math. Soc.} {\bf 26} (1994), no. 2, 140--146.
 

 

\bibitem{EliKer:1998a} S. Eliahou and M. Kervaire, Sumsets in vector spaces over finite fields. {\em J. Number Theory} {\bf 71} (1998), no. 1, 12--39.


 
\bibitem{Erd:1965a}  P. Erd\H{o}s, Extremal problems in number theory.  {\em Proc. Sympos. Pure Math., Vol. VIII} pp. 181--189 {\em  Amer. Math. Soc., Providence, R.I.}, 1965.
 




\bibitem{ErdGra:1980b} P. Erd\H{o}s and R. L. Graham, On bases with an exact order. {\em Acta Arith.} {\bf 37} (1980), 201--207.

\bibitem{ErdHei:1964a} P. Erd\H{o}s and H. Heilbronn, On the addition of residue classes mod $p$.  {\em Acta Arith.} {\bf 9} (1964), 149--159.


\bibitem{ErdNat:1987a} P. Erd\H{o}s and M. B. Nathanson, Problems and results on minimal bases in additive number theory. {\em Lecture Notes in Math.} {\bf 1240}, Springer, Berlin (1987), 87--96.

\bibitem{ErdTur:1941a} P. Erd\H{o}s and  P. Tur\'an, On a problem of Sidon in additive number theory and some related questions. {\em J. London Math. Soc.} {\bf 16} (1941), 212--215.



\bibitem{GalGre:2002a} L. Gallardo, G. Grekos,  L. Habsieger, F. Hennecart, B. Landreau, and A. Plagne, Restricted addition in $\mathbb{Z}/n\mathbb{Z}$ and an application to the Erd\H{o}s--Ginzburg--Ziv problem.  {\em J. London Math. Soc. (2)} {\bf 65} (2002), no. 3, 513--523.

 
\bibitem{GreRuz:2005a} B. Green and I. Ruzsa, sumfree sets in abelian groups.  {\em Israel J. Math.} {\bf 147} (2005), 157--188.
 
\bibitem{HamPla:2003a} Y. O. Hamidoune and A. Plagne, A new critical pair theorem applied to sumfree sets in abelian groups.  {\em Comment. Math. Helv.} {\bf 79} (2004), no. 1, 183--207.
 


 
\bibitem{Kar:2003a} Gy. K\'arolyi, On restricted set addition in abelian groups.  {\em Ann. Univ. Sci. Budapest, E\"otv\"os Sect. Math.} {\bf 46}, (2003) 47--54 (2004).
 
\bibitem{Kar:2004a} Gy. K\'arolyi, The Erd\H{o}s--Heilbronn problem in abelian groups.  {\em Israel J. Math.} {\bf 139} (2004) 349--359.
 
\bibitem{Kar:2005a} Gy. K\'arolyi, An inverse theorem for the restricted set addition in abelian groups.  {\em J. Algebra} {\bf 290} (2005), no. 2, 557--593.
 
 
 
\bibitem{Kem:1960a} J. H. B. Kemperman, On small sumsets in an abelian group.  {\em Acta Math.} {\bf 103} (1960), 63--88.
 

\bibitem{Kne:1953a}  M. Kneser, Absch\"atzungen der asymptotichen Dichte von Summenmengen. {\em Math. Z.} {\bf 58} (1953). 459--484.

 




\bibitem{LamThaPla:2017a} V. Lambert, T. H. L\^e, and A. Plagne, Additive bases in groups.  {\em Israel J. Math.}  {\bf 217} (2017), no. 1, 383--411.

\bibitem{Lev:2000a} V. F. Lev, Restricted set addition in groups. I. The classical setting.  {\em J. London Math. Soc. (2)} {\bf 62} (2000), no. 1, 27--40.



\bibitem{Nag:1948a} T. Nagell, L\o{}sning til oppgave nr. 2, 1943, s. 29.  {\em Nordisk Mat. Tidskr.} {\bf 30} (1948), 62--64.


\bibitem{Nag:1961a} T. Nagell, The Diophantine Equation $x^2+7=2^n$. {\em Ark. Mat.} {\bf 4} (1961), 185--187.


 

\bibitem{Nat:1974a} M. B. Nathanson, Minimal bases and maximal nonbases in additive number theory. {\em J. Number Theory} {\bf 6} (1974), 324--333.


 
\bibitem{Nat:2014a} M. B. Nathanson, Paul Erd\H{o}s and additive bases.  arXiv:1401.7598 [math.NT] (2014).





\bibitem{Nat:1996a} M. B. Nathanson, Additive number theory. Inverse problems and the geometry of sumsets.  Graduate Texts in Mathematics, {\bf 165}. {\em Springer--Verlag, New York},  1996. xiv+293 pp.
 
 

\bibitem{Pla:2002a} A. Plagne, Maximal $(k,l)$-free sets in $\mathbb{Z}/p\mathbb{Z}$ are arithmetic progressions. {\em Bull. Austral. Math. Soc.} {\bf 65} (2002), no. 3, 137--144.

 
\bibitem{Pla:2006a} A. Plagne, Optimally small sumsets in groups, I. The supersmall sumset property, the $\mu_G^{(k)}$ and the $\nu_G^{(k)}$ functions. {\em Unif. Distrib. Theory} {\bf 1} (2006), no. 1, 27--44.
 


\bibitem{Ram:2000a} S. Ramanujan, Question 464. {\em Journal of the Indian Mathematical Society} {\bf 5} (1913), 120.



 
\bibitem{Vos:1956a} A. G. Vosper, The critical pairs of subsets of a group of prime order.  {\em J. London Math. Soc.} {\bf 31} (1956), 200--205.
 
\bibitem{Vos:1956b} A. G. Vosper,  Addendum to ``The critical pairs of subsets of a group of prime order''. {\em J. London Math. Soc.} {\bf 31} (1956), 280--282.
 
\bibitem{WalStrWal:1972a} W. D. Wallis, A. P. Street, and J. S. Wallis, Combinatorics: room squares, sumfree sets, Hadamard matrices. Lecture Notes in Mathematics, {\bf 292}, {\em Springer--Verlag, Berlin-New York}, 1972. iv+508 pp.






















\end{thebibliography}
\end{document}